\documentclass[a4paper,12pt]{article}      % arXiv 

\usepackage{latexsym}

\usepackage{amsmath}
\usepackage{amsthm}
\usepackage{amssymb}
\usepackage{graphicx}
\usepackage{color} 

\title{
Note on Steepest Descent Algorithm for
\\
Quasi \Lnat-convex Function Minimization% 
\footnote{
This work was supported by JSPS KAKENHI Grant Numbers 
JP23K11001 and JP23K10995.
}%
}%

\author{Kazuo Murota%
\footnote{
The Institute of Statistical Mathematics,
Tokyo 190-8562, Japan; 
and
Faculty of Economics and Business Administration,
Tokyo Metropolitan University, 
Tokyo 192-0397, Japan,
\texttt{murota@tmu.ac.jp}
}
 \and
Akiyoshi Shioura%
\footnote{
Department of Industrial Engineering and Economics,
Tokyo Institute of Technology,
Tokyo 152-8550, Japan
  \texttt{shioura.a.aa@m.titech.ac.jp}}
}

\date{July 2023} %% For arXiv

\newcommand{\RED}[1]{{\color{red}#1}}
\newcommand{\BLU}[1]{{\color{blue}#1}}
 \renewcommand{\RED}[1]{{#1}} %% For arXiv
 \renewcommand{\BLU}[1]{{#1}} %% For arXiv

\newcommand{\OMIT}[1]{{\bf [OMIT:} #1 \ {\bf --- end OMIT] }}  %%% For work
   \renewcommand{\OMIT}[1]{}            %%% For FINAL
\usepackage{geometry}
\newcommand{\suppp}{{\rm supp}\sp{+}}
\newcommand{\suppm}{{\rm supp}\sp{-}}

\newcommand{\RR}{\mathbb{R}}

\newcommand{\Z}{\mathbb{Z}}

\newcommand{\dom}{{\rm dom\,}}
\newcommand{\unitvec}[1]{e_{#1}}

\newcommand{\Rinf}{\RR \cup \{ +\infty \}}

\newcommand{\Lnat}{{L$^{\natural}$}}

\newcommand{\1}{\mathbf{1}}

\newtheorem{theorem}{Theorem}[section]
\newtheorem{lemma}[theorem]{Lemma}

\numberwithin{equation}{section}

\begin{document}

\maketitle

\begin{abstract}
We define a class of discrete quasi convex functions,
called semi-strictly quasi \Lnat-convex functions,
and show that the steepest descent algorithm
for \Lnat-convex function minimization
also works for this class of quasi convex functions.
The analysis of the exact number of iterations
is also extended, revealing the so-called geodesic property
of the steepest descent algorithm when applied to 
semi-strictly quasi \Lnat-convex functions.
\end{abstract}

\section{Results}

We define a class of discrete quasi convex functions,
called semi-strictly quasi \Lnat-convex functions,
and extend the results of Murota--Shioura \cite{MS14exbndLmin}
on the steepest descent algorithm
designed for \Lnat-convex function minimization.
Specifically, we show that the steepest descent algorithm
also works for this class of quasi convex functions
with the so-called geodesic property.
The exact number of iterations required by the steepest descent algorithm 
applied to semi-strictly quasi \Lnat-convex functions
is measured by a variant of the L$_\infty$-distance
between the initial point and the nearest minimizer,
similarly for (ordinary) \Lnat-convex function minimization.

\medskip

 We say that a function $g: \Z^n \to \Rinf$ is 
{\it semi-strictly quasi \Lnat-convex} 
(s.s.~quasi \Lnat-convex) 
if for every $p, q \in \dom g$ and every nonnegative $\lambda \in \Z_+$,
\RED{%
both (I) and (II) below are satisfied.
\\
(I)
One of the following three conditions holds:
\begin{align}
&  g(p) > g((p - \lambda \1) \vee  q),
 \label{eqn:def-Lnat1P}
\\
&   g(q) > g(p  \wedge (q + \lambda \1)),
 \label{eqn:def-Lnat2P}
\\
&  g(p) =g((p - \lambda \1) \vee q) \mbox{ and }
 g(q) = g(p  \wedge (q + \lambda \1)).
 \label{eqn:def-Lnat3P}
\end{align}
(II)
One of the following three conditions holds:
\begin{align}
&  g(q) > g((p - \lambda \1) \vee  q),
 \label{eqn:def-Lnat1Q}
\\
&   g(p) > g(p  \wedge (q + \lambda \1)),
 \label{eqn:def-Lnat2Q}
\\
&  g(q) =g((p - \lambda \1) \vee q) \mbox{ and }
 g(p) = g(p  \wedge (q + \lambda \1)).
 \label{eqn:def-Lnat3Q}
\end{align}
}%red
Here $\dom g= \{p \in \Z^n \mid g(p) < + \infty\}$,
$\1 = (1,1,\ldots, 1)$,
and for $p, q \in \Z^n$ the vectors $p \wedge q$
and $p \vee q$ denote, respectively, the vectors
of component-wise minimum and maximum of $p$ and $q$.
\RED{%
It is noted that, 
while quasi L-convexity was investigated
in Murota--Shioura \cite{MS03quasi} 
(see also Murota \cite[Section 7.11]{Mdcasiam}),
there is no formulation of the concept of 
quasi \Lnat-convexity in the literature. 
The above definition of s.s.~quasi \Lnat-convex functions
is an ordinal version of the translation-submodularity
(see \cite[Section 7.1]{Mdcasiam})
of \Lnat-convex functions. 
}%

 We consider minimization of 
a semi-strictly quasi \Lnat-convex function 
$g:\Z^n \to \Rinf$ 
with $\arg\min g \neq \emptyset$.
 Let $N = \{1, 2, \ldots, n\}$ and 
denote by 
\RED{%
$\unitvec{X} \in \{0,1\}^n$ 
}%
the characteristic
vector of $X \subseteq N$, i.e.,
\RED{%
$\unitvec{X}(i)=1$
}%
if $i \in X$ and 
\RED{%
$\unitvec{X}(i)=0$
}%
if $i \in N \setminus X$.
The local optimality condition for s.s.~quasi \Lnat-convex functions 
is exactly the same as that for \Lnat-convex functions.

\begin{theorem}
\label{thm:opt-cond}
\RED{%
Let $g:\Z^n \to \Rinf$ be a semi-strictly quasi \Lnat-convex function.
}%
A vector $p \in \dom g$ is a minimizer of $g$ if and only if
$g(p+ \sigma \unitvec{X}) \ge g(p)$ holds for every $\sigma \in \{+1, -1\}$
and $X \subseteq N$.
\end{theorem}

\begin{proof}
The proof is given in Section \ref{SCprfTHopt-cond}.
\end{proof}

The following is the basic form of 
the steepest descent algorithm
for \Lnat-convex function minimization
(see \cite[Section 10.3.1]{Mdcasiam}).

\medskip

\noindent
{\bf Algorithm} {\sc SteepestDescent}\\
Step 0: 
 Set $p:=p^\circ$.
\\
Step 1:
Find $\sigma \in \{+1, -1\}$
and $X \subseteq N$ that minimize $g(p+ \sigma \unitvec{X})$.
\\
Step 2:
 If $g(p + \sigma \unitvec{X}) = g(p)$, then output $p$ and stop.
\\
Step 3:
 Set $p:=p + \sigma \unitvec{X}$ and go to Step 1.

\medskip

 In this technical note, 
we first show that the output of Algorithm {\sc SteepestDescent}
is a minimizer of $g$
\RED{%
under a weaker condition that
$g$ is an s.s.~quasi \Lnat-convex function,
and that the exact number of iterations
revealed by Murota--Shioura \cite{MS14exbndLmin}
for \Lnat-convex functions
remains the same.
}%
The number of iterations 
is measured by the ``distance'' between the initial vector and a minimizer of $g$.
 For a vector $q \in \Z^n$, denote
\begin{equation}
\label{normpmdef}
\|q\|_\infty^+ = \max_{i \in N}\, \max(0, q(i)),
\qquad
\|q\|_\infty^- = \max_{i \in N}\, \max(0, -q(i)).
\end{equation}
 Note that
\[
 \|q\|_\infty
=  \max(\|q\|_\infty^+ , \|q\|_\infty^- )
\]
holds, and $\|q\|_\infty^+ + \|q\|_\infty^-$ serves as
a norm of $q$ (satisfying the axioms of norms).
 Accordingly, the value $\|p^* -p \|_\infty^+ + \|p^* - p \|_\infty^-$
represents a distance between two vectors $p^*$  and $p$.
For $p \in \Z^n$, we define
\begin{equation} \label{mupdef}
 \mu(p)
 =   \min\{\|p^* -p \|_\infty^+
  + \|p^* - p \|_\infty^- \mid p^* \in \arg\min g\},
\end{equation}
which measures the distance between the vector $p$ and the set of minimizers of $g$.

 It is easy to see that $\mu(p)$ decreases by at most one
if $p$ is updated by adding or subtracting a 0-1 vector,
i.e., $\mu(p + \sigma \unitvec{X}) \geq \mu(p)-1$ for $\sigma \in \{+1, -1\}$ and
$X\subseteq N$.
 This implies that $\mu(p)+1$ is a lower bound for the number of
 iterations in {\sc SteepestDescent}.
 This is also  an upper bound as follows.

\begin{theorem}
\label{thm:greedy-iter}
 The algorithm {\sc SteepestDescent},
\RED{%
when applied to a semi-strictly quasi \Lnat-convex function $g$,
}%
outputs a minimizer of $g$ and 
terminates exactly in ${\mu}(p^\circ)+1$ iterations.
\end{theorem}

\begin{proof}
The proof is given in Section \ref{SCprfTHgreedy-iter}.
\end{proof}

Theorem \ref{thm:greedy-iter} shows 
\BLU{%
the so-called geodesic property
}%
that the trajectory of a vector $p$
generated by the steepest descent algorithm 
is the ``shortest'' path 
between the initial vector
and a minimizer of $g$.

The following variant of the steepest descent algorithm,
where the vector $p$ is always incremented,
\RED{%
has been proposed in 
\cite{MS14exbndLmin}
for \Lnat-convex functions.
We next show that this variant also works for 
s.s.~quasi \Lnat-convex functions.
}%

\medskip

\noindent
{\bf Algorithm} {\sc SteepestDescentUp}%
\\
Step 0:  Set $p:=p^\circ$, where
$p^\circ \in \Z^n$ is a lower bound of some $p^* \in \arg\min g$.\\
Step 1:
Find $X \subseteq N$ that minimizes $g(p+ \unitvec{X})$.
\\
Step 2:
 If $g(p + \unitvec{X}) = g(p)$, then output $p$ and stop.
\\
Step 3:
 Set $p:=p + \unitvec{X}$ and go to Step 1.

\medskip

 For the analysis of {\sc SteepestDescentUp}, we define
\begin{equation} \label{hatmupdef}
 \hat\mu(p) = \min\{\|p^* - p\|_\infty \mid p^* \in \arg\min g,\ p^* \geq p \}
\quad (p \in \Z^n).
\end{equation}

\begin{theorem}
\label{thm:ascend-iter}
\RED{%
Let $g:\Z^n \to \Rinf$ be a semi-strictly quasi \Lnat-convex function, and 
}%
%% Let $g:\Z^n \to \Rinf$ be an \Lnat-convex function with $\arg\min g \neq \emptyset$.
suppose that the initial vector $p^\circ \in \dom g$ in
the algorithm {\sc SteepestDescentUp} is a lower bound of
some minimizer of $g$.
 Then, the algorithm outputs a minimizer of $g$ and
terminates exactly in $\hat{\mu}(p^\circ)+1$ iterations.
\end{theorem}
\begin{proof}
The proof is given in Section \ref{SCprfTHascend-iter}.
\end{proof}

 Similarly to {\sc SteepestDescentUp},
we can consider an algorithm {\sc SteepestDescentDown}
\RED{%
\cite{MS14exbndLmin},
}%
where the vector is decreased by a vector $\unitvec{X} \in \{0,1\}^n$
that minimizes $g(p-\unitvec{X})$.

\begin{theorem}
\label{thm:descend-iter}
\RED{%
Let $g:\Z^n \to \Rinf$ be a semi-strictly quasi \Lnat-convex function, and 
}%
suppose that the initial vector $p^\circ {\in} \dom g$ in
the algorithm {\sc SteepestDescentDown} is an upper bound of
some minimizer of $g$.
 Then, the algorithm outputs a minimizer of $g$ and
terminates exactly in $\check{\mu}(p^\circ)+1$
iterations, where
\begin{equation} \label{checkmupdef}
 \check\mu(p) = \min\{\|p^* - p\|_\infty \mid p^* \in \arg\min g,\ p^* \leq p \}
\quad (p \in \Z^n).
\end{equation}
\end{theorem}
\begin{proof}
It follows from the definition of an s.s.~quasi \Lnat-convex function
that a function $g$ is s.s.~quasi \Lnat-convex
 if and only if the function $h(p):=g(-p)$ is s.s.~quasi \Lnat-convex.
In addition, 
the algorithm {\sc SteepestDescentDown} 
applied to $g$ with an initial vector $p^\circ$ 
behaves `isomorphically' to 
the algorithm {\sc SteepestDescentUp} 
applied to $h$ with initial vector $-p^\circ$. 
Therefore, Theorem \ref{thm:descend-iter}
follows from Theorem \ref{thm:ascend-iter}.
\end{proof}

%%\newpage

\section{Proofs}

 In this section, we prove 
Theorems \ref{thm:opt-cond},
\ref{thm:greedy-iter},
and \ref{thm:ascend-iter}.
The key facts used in the proofs
\RED{%
are the following properties
}%
 of s.s.~quasi \Lnat-convex functions.
 For $p \in \Z^n$, we denote $\suppp(p) = \{i \in N \mid p(i)>0\}$.

\begin{lemma}
\label{lem:Lnat-exchange}
Let $g:\Z^n \to \Rinf$ be 
\RED{%
a semi-strictly quasi \Lnat-convex function.
}%
 For every $p, q \in \dom g$ with $\suppp(p-q) \neq \emptyset$,
one of the following three conditions holds
with $Y = \arg\max_{i \in N}\{p(i) - q(i)\}:$
\begin{align}
& 
\label{eqn:lem:Lnat-exchange-1}
 g(p) > g(p - \unitvec{Y}),
\\
& 
\label{eqn:lem:Lnat-exchange-2}
 g(q) > g(q + \unitvec{Y}),
\\
& 
\label{eqn:lem:Lnat-exchange-3}
 g(p) = g(p - \unitvec{Y}) \mbox{ and } g(q) = g(q + \unitvec{Y}).
\end{align}
\end{lemma}

\begin{proof}
 Let $p, q \in \dom g$, and suppose that  $\suppp(p-q) \neq \emptyset$.
\RED{%
By \eqref{eqn:def-Lnat1Q}--\eqref{eqn:def-Lnat3Q} 
in the definition of an s.s.~quasi \Lnat-convex function,
}%
one of the following three conditions holds for every nonnegative $\lambda \in \Z_+$:
\begin{align*}
&  g(q) > g((p - \lambda \1) \vee  q),
%% \label{eqn:def-Lnat1Qpr}
\\
&   g(p) > g(p  \wedge (q + \lambda \1)),
%% \label{eqn:def-Lnat2Qpr}
\\
&  g(q) =g((p - \lambda \1) \vee q) \mbox{ and }
 g(p) = g(p  \wedge (q + \lambda \1)).
%% \label{eqn:def-Lnat3Qpr}
\end{align*}
Putting $\lambda = \max_{i \in N}\{p(i)-q(i)\} - 1$,
we have 
\[
 \lambda \geq 0, \quad
(p - \lambda \1) \vee q = q + \unitvec{Y},
\quad
p \wedge (q+ \lambda \1) = p - \unitvec{Y}.
\]
Substituting these expressions into the above, we obtain
\begin{align*}
&  g(q) > g((p - \lambda \1) \vee  q) = g(q + \unitvec{Y}),
\\
&   g(p) > g(p  \wedge (q + \lambda \1)) = g(p - \unitvec{Y}),
\\
&  g(q) =g(q + \unitvec{Y}) \mbox{ and }
 g(p) = g(p - \unitvec{Y}) .
\end{align*}
Thus we obtain \eqref{eqn:lem:Lnat-exchange-1}--\eqref{eqn:lem:Lnat-exchange-3}.
\end{proof}

\begin{lemma}  \label{lem:Lnat-subm}
Let $g:\Z^n \to \Rinf$ be a semi-strictly quasi \Lnat-convex function.
 For every $p, q \in \dom g$,
one of the following three conditions holds
with $Z = \suppp(p-q):$
\begin{align}
& 
\label{eqn:lem:Lnat-subm-1}
 g(p) > g((p \vee q) - \unitvec{Z}),
\\
& 
\label{eqn:lem:Lnat-subm-2}
 g(q) > g((p \wedge q) + \unitvec{Z}),
\\
& 
\label{eqn:lem:Lnat-subm-3}
 g(p) = g((p \vee q) - \unitvec{Z}) \mbox{ and } g(q) = g((p \wedge q) + \unitvec{Z}).
\end{align}
\end{lemma}

\begin{proof}
By \eqref{eqn:def-Lnat1P}--\eqref{eqn:def-Lnat3P} 
in the definition of an s.s.~quasi \Lnat-convex function,
one of the following three conditions holds for every nonnegative $\lambda \in \Z_+$:
\begin{align*}
&  g(p) > g((p - \lambda \1) \vee  q),
%% \label{eqn:def-Lnat1Ppr}
\\
&   g(q) > g(p  \wedge (q + \lambda \1)),
%% \label{eqn:def-Lnat2Ppr}
\\
&  g(p) =g((p - \lambda \1) \vee q) \mbox{ and }
 g(q) = g(p  \wedge (q + \lambda \1)).
%% \label{eqn:def-Lnat3Ppr}
\end{align*}
Putting $\lambda = 1$,
we have 
\[
(p - \lambda \1) \vee q = (p \vee q) - \unitvec{Z},
\quad
p \wedge (q+ \lambda \1) = (p \wedge q) + \unitvec{Z}.
\]
Substituting these expressions into the above, we obtain
\eqref{eqn:lem:Lnat-subm-1}--\eqref{eqn:lem:Lnat-subm-3}.
\end{proof}

\subsection{Proof of Theorem \ref{thm:opt-cond}}
\label{SCprfTHopt-cond}

 We prove the ``if'' part only since 
the ``only if'' part of Theorem \ref{thm:opt-cond} is obvious.
\RED{%
Namely, 
}%
we show that if $p \in \dom g$ is not a minimizer of 
\RED{%
an s.s.~quasi \Lnat-convex function $g$,
}%
then there exist some
 $\sigma \in \{+1, -1\}$ and $X \subseteq N$
such that $g(p+ \sigma \unitvec{X}) < g(p)$.

 Let $q^* \in \dom g$ be a minimizer of $g$
that minimizes 
\BLU{%
$\|q^* - p\|_\infty^+ + \|q^* - p\|_\infty^-$.
}%
\BLU{%
Since $p$ is not a minimizer,
we have $q^* \ne p$, that is,
$\suppp(p - q^*) \ne \emptyset$
or
$\suppp(q^* - p) \ne \emptyset$.%
}%%

\RED{%
First we consider the case of $\suppp(p - q^*) \ne \emptyset$.
}%
 By Lemma \ref{lem:Lnat-exchange}
\RED{%
for $(p, q^*)$,
}%
 we have one of the following three conditions  
with $Y = \arg\max_{i \in N}\{p(i) - q^*(i)\}$:
\begin{align}
& 
\label{eqn:thm:opt-cond-1}
 g(p) > g(p - \unitvec{Y}),
\\
& 
\label{eqn:thm:opt-cond-2}
 g(q^*) > g(q^* + \unitvec{Y}),
\\
& 
\label{eqn:thm:opt-cond-3}
 g(p) = g(p - \unitvec{Y}) \mbox{ and } g(q^*) = g(q^* + \unitvec{Y}).
\end{align}
The case \eqref{eqn:thm:opt-cond-2} is excluded since $q^*$
is a minimizer of $g$.
 By the choice of $q^*$,
the vector $q^* + \unitvec{Y}$ is not a minimizer of $g$
since it satisfies 
\BLU{%
$\|(q^* + \unitvec{Y}) - p\|_\infty^+ + \|(q^* + \unitvec{Y}) - p\|_\infty^- 
  <\|q^* - p\|_\infty^+ + \|q^* - p\|_\infty^-$.
}%
 Hence, the case \eqref{eqn:thm:opt-cond-3} is also excluded.
 This shows that the condition \eqref{eqn:thm:opt-cond-1}
is satisfied,
i.e., 
$g(p+ \sigma \unitvec{X}) < g(p)$ holds with $\sigma = -1$ and $X = Y$.

\RED{%
The second case of $\suppp(q^* -p) \ne \emptyset$
can be treated symmetrically as follows.
}%
 By Lemma \ref{lem:Lnat-exchange}
\RED{%
for $(q^*,p)$,
}%
we have one of the following three conditions  
with $Y = \arg\max_{i \in N}\{ q^*(i) - p(i)\}$:
\begin{align}
& 
\label{eqn:thm:opt-cond-1Q}
 g(p) > g(p + \unitvec{Y}),
\\
& 
\label{eqn:thm:opt-cond-2Q}
 g(q^*) > g(q^* - \unitvec{Y}),
\\
& 
\label{eqn:thm:opt-cond-3Q}
 g(p) = g(p + \unitvec{Y}) \mbox{ and } g(q^*) = g(q^* - \unitvec{Y}).
\end{align}
The case \eqref{eqn:thm:opt-cond-2Q} is excluded since $q^*$
is a minimizer of $g$.
 By the choice of $q^*$,
the vector $q^* - \unitvec{Y}$ is not a minimizer of $g$
since it satisfies 
\BLU{%
$\|(q^* - \unitvec{Y}) - p\|_\infty^+  + \|(q^* - \unitvec{Y}) - p\|_\infty^- 
< \|q^* - p\|_\infty^+ + \|q^* - p\|_\infty^-$.
}%
 Hence, the case \eqref{eqn:thm:opt-cond-3Q} is also excluded.
 This shows that the condition \eqref{eqn:thm:opt-cond-1Q}
is satisfied,
i.e., 
$g(p + \sigma \unitvec{X}) < g(p)$ holds with $\sigma = +1$ and $X = Y$.

\subsection{Proof of Theorem \ref{thm:greedy-iter}}
\label{SCprfTHgreedy-iter}

 The bound $\mu(p^\circ) + 1$ for the number of iterations in
algorithm {\sc SteepestDescent}
\RED{%
applied to an s.s.~quasi \Lnat-convex function
}%
can be obtained by repeated application of the following lemma.

\begin{lemma}
\label{lem:Lnat-greedy}
 Let $p \in \Z^n$ be a vector with $\mu(p) > 0$.
 Suppose that $\sigma \in \{+1, -1\}$
and $X \subseteq N$ minimize the value $g(p+ \sigma \unitvec{X})$.
 Then, 
$\mu(p+  \sigma\unitvec{X})= \mu(p)-1$.
\end{lemma}

\RED{%
To prove Lemma \ref{lem:Lnat-greedy},
we distinguish two cases,
depending on $\sigma = + 1$ or $\sigma = - 1$.
First we deal with the case of $\sigma = + 1$
while the other case of $\sigma = - 1$ is treated later
(by a symmetric argument)
}%%

\RED{%
 {\bf [Step 1 ($\sigma = + 1$)]} \quad
}%
 We first show the inequality
$\mu(p + \unitvec{X}) \geq \mu(p) -1$.
 For every $d \in \Z^n$ and $Y \subseteq N$, we have
\begin{equation*}
    \|d - \unitvec{Y} \|^+_\infty \geq  \|d  \|^+_\infty -1,\qquad
    \| d - \unitvec{Y} \|^-_\infty \geq  \| d  \|^-_\infty.
\end{equation*}
 Hence, it holds that
\begin{eqnarray*}
\mu(p + \unitvec{X})
& = &
\min\{ \|q - (p + \unitvec{X}) \|^+_\infty
+ \|q - (p + \unitvec{X}) \|^-_\infty
\mid q \in \arg\min g \}
\nonumber \\
& \geq &
\min\{ \| q - p \|^+_\infty +  \| q - p \|^-_\infty
 \mid q \in \arg\min g \} - 1
\nonumber \\
& = &  \mu(p) -1.
\end{eqnarray*}

\RED{%
In the following (i.e., in Steps 2 to 4),
}%
we prove the reverse inequality
\begin{equation}
\label{eqn:rev-ineq} 
\mu(p + \unitvec{X}) \leq \mu(p) -1.
\end{equation} 
The outline of the proof is as follows.
 We denote
\begin{align}
  S & = \{ q \in \arg\min g \mid
 \|q - p \|^+_\infty +  \| q - p \|^-_\infty = \mu(p) \},
\\
 \xi & = \max\{\|q - p\|_\infty^+ \mid q \in S \}.
\end{align}
 Let  $q^*$ be a vector in $S$
with $\|q^* - p \|^+_\infty = \xi$,
and assume that $q^*$ is a minimal vector
among all such vectors.
\RED{%
In Step 2,
}%
we show that 
\RED{%
\begin{equation}
\label{eqn:claim1} 
\xi > 0,  \quad  \mbox{i.e.,} \quad  \|q^* - p \|^+_\infty > 0.
\end{equation}
}%
 Note that this condition is equivalent to
$\suppp(q^* - p)\neq \emptyset$.
 Using this, we then prove,
\RED{%
in Step 3,
}%
 that
\begin{equation}
\label{eqn:claim2}
\arg\max_{i \in N}\{q^*(i) - p(i)\} \subseteq X. 
\end{equation}
 By using (\ref{eqn:claim1}) and (\ref{eqn:claim2}), 
we derive the inequality (\ref{eqn:rev-ineq})
\RED{%
in Step 4.
}%

\RED{%
 {\bf [Step 2  ($\sigma = + 1$): proof of (\ref{eqn:claim1})]} \quad
}%
 Assume, to the contrary, that  
\RED{%
$(\xi=) \ \|q^* - p \|^+_\infty = 0$, 
}%
i.e., $q^* \leq  p$ holds. 
 This assumption implies $\|q^* - p \|^-_\infty > 0$
since $\mu(p) > 0$.
\BLU{%
By Lemma \ref{lem:Lnat-subm}
for $(p+ \unitvec{X},q^*)$,
}%%
one of the following three conditions holds:
\BLU{%
\begin{align}
& 
 g(p + \unitvec{X}) > g(((p + \unitvec{X}) \vee q^*) - \unitvec{Z}),
 \label{eqn:lem-greedy7-1-1}
\\
& 
 g(q^*) > g(((p + \unitvec{X}) \wedge q^*) + \unitvec{Z}),
 \label{eqn:lem-greedy7-1-2}
\\
& 
 g(p + \unitvec{X}) = g(((p + \unitvec{X}) \vee q^*) - \unitvec{Z}) 
\mbox{ and } g(q^*) = g(((p + \unitvec{X}) \wedge q^*) + \unitvec{Z}),
 \label{eqn:lem-greedy7-1-3}
\end{align}
where $Z= \suppp((p + \unitvec{X})- q^*)$.
}%%
 Let 
\BLU{%
$Y = \{i \in N \mid p(i) - q^*(i) = 0\}$, 
}%
which may be the empty set.
Since 
\BLU{%
$p \geq q^*$ by the assumption
and $Z = (N \setminus Y) \cup X$,
}%% 
we have
\BLU{%
\begin{align*}
& 
((p + \unitvec{X}) \vee q^*) - \unitvec{Z}
 = (p + \unitvec{X})  - \unitvec{Z}
 = p - \unitvec{N \setminus (X\cup Y)},
\\ &
((p + \unitvec{X}) \wedge q^*) + \unitvec{Z}=
 q^* + \unitvec{Z}=
q^* + \unitvec{(N \setminus Y) \cup X}.
\end{align*}
}%color
With these equations, the three conditions 
 (\ref{eqn:lem-greedy7-1-1})--(\ref{eqn:lem-greedy7-1-3})
can be rewritten as follows:
\begin{align}
& 
  g(p+ \unitvec{X}) >  g(p - \unitvec{N \setminus (X \cup Y)}),
 \label{eqn:lem-greedy7-2-1}
\\
& 
  g(q^*)  >   g(q^* + \unitvec{(N \setminus Y) \cup X}),
 \label{eqn:lem-greedy7-2-2}
\\
& 
  g(p+ \unitvec{X}) =  g(p - \unitvec{N \setminus (X \cup Y)})
\mbox{ and }
  g(q^*)  =   g(q^* + \unitvec{(N \setminus Y) \cup X}).
 \label{eqn:lem-greedy7-2-3}
\end{align} 
\RED{%
By the choice of $(\sigma,X)$,
where 
$\sigma =+1$ in our first case, 
}% 
we have
$g(p+ \unitvec{X}) \leq
g(p - \unitvec{N \setminus (X \cup Y)})$.
 From this and (\ref{eqn:lem-greedy7-2-1})--(\ref{eqn:lem-greedy7-2-3})
 follows that
$ g(q^*)  \geq
g(q^* + \unitvec{(N \setminus Y) \cup X})$,
implying that
$q^* + \unitvec{(N \setminus Y) \cup X} \in \arg\min g$.
 By $q^* \leq p$ and the definition of $Y$, we have
\begin{equation}
 \label{eqn:lem-greedy7-11} 
\|(q^* + \unitvec{(N \setminus Y) \cup X}) - p\|_\infty^-
 =  
\max_{i \in N \setminus Y}\{p(i) - (q^*(i)+1)\}
 =  
\|q^* - p\|_\infty^- -1
\end{equation}
since $\|q^* - p\|_\infty^- > 0$.
 We also have
\begin{align}
\|(q^* + \unitvec{(N \setminus Y) \cup X}) - p\|_\infty^+
& \leq  1\  =\  \|q^* - p\|_\infty^+ +1 
 \label{eqn:lem-greedy7-3} 
\end{align}
since $\|q^* - p\|_\infty^+ = 0$.
 From (\ref{eqn:lem-greedy7-11}) and (\ref{eqn:lem-greedy7-3}) follows that
\begin{align}
\mu(p)
& \leq 
\|(q^* + \unitvec{(N \setminus Y) \cup X})  - p \|^+_\infty
 + \|(q^* + \unitvec{(N \setminus Y) \cup X}) - p\|_\infty^-
\notag\\
&\leq   \|q^*  - p \|^+_\infty + \|q^* - p\|_\infty^-
 =  \mu(p),
 \label{eqn:lem-greedy7-12} 
\end{align}
where the first inequality is by 
the definition of $\mu(p)$.
 Hence, the inequality (\ref{eqn:lem-greedy7-3}) 
and the first inequality in (\ref{eqn:lem-greedy7-12}) must hold with equality,
i.e., we have $q^* + \unitvec{(N \setminus Y) \cup X} \in S$ 
and 
\begin{equation*}
 \|(q^* + \unitvec{(N \setminus Y) \cup X}) - p\|_\infty^+
 =  \|q^* - p\|_\infty^+ +1 
>  \|q^* - p\|_\infty^+ = \xi. 
\end{equation*}
 This, however, is a contradiction to the definition of $\xi$.
 Hence, (\ref{eqn:claim1}) holds.

\RED{%
 {\bf [Step 3  ($\sigma = + 1$): proof of (\ref{eqn:claim2})]} \quad
}%
 We denote
\begin{equation*}
A = \arg\max_{i \in N}\{q^*(i) - p(i)\}.
\end{equation*}
 Then, (\ref{eqn:claim2}) is simply rewritten as $A \subseteq X$.

 Assume, to the contrary, that $A \setminus X \neq \emptyset$ holds.
\RED{%
 We  will show 
}%
that
$q^* - \unitvec{A \setminus X} \in \arg\min g$.
 By (\ref{eqn:claim1})
\RED{%
established in Step 2,
}%
it holds that  $\xi=\|q^* - p\|_\infty^+   > 0$.
 Therefore, we have $A \subseteq \suppp(q^* - p)$, from which follows that 
\begin{equation*}
 \suppp(q^* - (p + \unitvec{X})) \supseteq A \setminus X
\neq \emptyset.
\end{equation*}
 Since $A \setminus X \neq \emptyset$, we also have
\begin{equation*}
 \arg\max_{i \in N}\{q^*(i) - (p + \unitvec{X})(i)\} = A \setminus X.
\end{equation*}
 Hence, Lemma \ref{lem:Lnat-exchange} 
\RED{%
for $(q^*, p+ \unitvec{X})$,
together with the relation
$(p + \unitvec{X}) + \unitvec{A \setminus X} = p + \unitvec{X \cup A}$,
implies that
}%
one of the following three conditions holds:
\begin{align}
&
   g(q^*)  > g(q^* - \unitvec{A \setminus X}),
 \label{eqn:lem-greedy2-1}
\\
&
 g(p+ \unitvec{X}) > g(p + \unitvec{X \cup A}),
 \label{eqn:lem-greedy2-2}
\\
&
   g(q^*)  = g(q^* - \unitvec{A \setminus X})
\mbox{ and }
 g(p+ \unitvec{X}) = g(p + \unitvec{X \cup A}).
 \label{eqn:lem-greedy2-3}
\end{align}
\RED{%
 By the choice of $(\sigma,X)$,
where 
$\sigma =+1$ in our first case, 
}% 
 we have
$g(p+ \unitvec{X}) \leq  g(p + \unitvec{X \cup A})$,
which, together with (\ref{eqn:lem-greedy2-1})--(\ref{eqn:lem-greedy2-3}), 
implies that
 $g(q^*) \geq g(q^* - \unitvec{A \setminus X})$, i.e.,
$q^* - \unitvec{A \setminus X} \in \arg\min g$.

 Since $A  \setminus X \subseteq A \subseteq \suppp(q^* - p)$, we have
\begin{align}
 \|(q^*- \unitvec{A \setminus X}) - p\|_\infty^+
& \leq  \|q^* - p\|_\infty^+ = \xi,
 \label{eqn:lem-greedy7-5} 
\\
\|(q^*- \unitvec{A \setminus X}) - p\|_\infty^-
& = \|q^* - p\|_\infty^-,
%%\notag
% \label{eqn:lem-greedy7-8} 
\end{align}
 from which follows that
\begin{align}
 \mu(p)
& \leq \|(q^*- \unitvec{A \setminus X}) - p\|_\infty^+
 +
 \|(q^*- \unitvec{A \setminus X}) - p\|_\infty^- \notag\\
& \leq \|q^* - p\|_\infty^+
 + \|q^* - p\|_\infty^- = \mu(p),
  \label{eqn:lem-greedy7-6} 
\end{align}
where the first inequality is by the definition of $\mu(p)$. 
 Hence, the inequality (\ref{eqn:lem-greedy7-5})
and the first inequality in (\ref{eqn:lem-greedy7-6}) 
must hold with equality.
 Hence, the vector $q^* - \unitvec{A \setminus X}$
belongs to $S$ with
$\|(q^* - \unitvec{A \setminus X}) - p\|_\infty^+  = \xi$,
a contradiction to the minimality of $q^*$.

\RED{%
 {\bf [Step 4  ($\sigma = + 1$): proof of (\ref{eqn:rev-ineq})]} \quad
}%
 To derive the inequality (\ref{eqn:rev-ineq})
\RED{%
from (\ref{eqn:claim1}) and (\ref{eqn:claim2}),
 we distinguish two cases:
$\min_{i \in N}\{q^*(i) - p(i)\} > 0$ 
or
$\min_{i \in N}\{q^*(i) - p(i)\} \leq 0$.%
}%

We first consider the case 
with $\min_{i \in N}\{q^*(i) - p(i)\} > 0$.
%%or $B \cap X = \emptyset$ holds.
 Since $q^*(i) > p(i)$ for all $i \in N$, it holds that
$q^* \geq p + \unitvec{X}$.
 Therefore, we have
\begin{equation*}
  \|q^* - (p + \unitvec{X}) \|^-_\infty
=0= \|q^* - p \|^-_\infty.
\end{equation*}
 By (\ref{eqn:claim2})
\RED{%
established in Step 3,
}%,
it holds that
\begin{equation*}
  \|q^* - (p + \unitvec{X}) \|^+_\infty
= \|q^* - p \|^+_\infty - 1.
\end{equation*}
 Therefore, it follows that
\begin{align*}
\mu(p+ \unitvec{X}) 
& \leq 
  \|q^* - (p + \unitvec{X}) \|^+_\infty
+ \|q^* - (p + \unitvec{X}) \|^-_\infty
\nonumber
\\
& = 
  (\|q^* - p  \|^+_\infty -1) + \|q^* - p  \|^-_\infty 
 =  \mu(p) -1.
\end{align*}

We next consider the remaining case where
 $\min_{i \in N}\{q^*(i) - p(i)\} \leq 0$.
%and  $B \cap X \neq \emptyset$.
 We denote
\begin{equation*}
 B = \arg\min_{i \in N}\{q^*(i) - p(i)\}.
\end{equation*}
\RED{%
 We  will show 
}%
that
$q^* + \unitvec{B \cap X} \in \arg\min g$ holds.
 If $B \cap X = \emptyset$, then
$q^* + \unitvec{B \cap X} = q^*\in \arg\min g$.
 Hence, we assume
$B \cap X \neq \emptyset$.
 Since
\begin{equation*}
\max_{i \in N}\{p(i) - q^*(i)\} \geq 0,
\qquad  B = \arg\max_{i \in N}\{p(i) - q^*(i)\},
\end{equation*}
it holds that
\begin{equation}
\suppp((p+ \unitvec{X})- q^*) \supseteq B \cap X \neq \emptyset, 
\quad
 \arg\max_{i \in N}\{(p+ \unitvec{X})(i)- q^*(i)\} = B \cap X.
\end{equation}
 It follows from Lemma \ref{lem:Lnat-exchange} 
\RED{%
for $(p+ \unitvec{X}, q^*)$
}%
and the equation
$(p + \unitvec{X}) - \unitvec{B \cap X} = p + \unitvec{X \setminus B}$
that one of the following three conditions holds:
\begin{align}
&
   g(p+ \unitvec{X})    >  g(p + \unitvec{X \setminus B}),
 \label{eqn:lem-greedy3-1}
\\
&
g(q^*)  > g(q^* + \unitvec{B \cap X}),
 \label{eqn:lem-greedy3-2}
\\
&
   g(p+ \unitvec{X})    =  g(p + \unitvec{X \setminus B}) 
\mbox{ and }
g(q^*)  = g(q^* + \unitvec{B \cap X}).
 \label{eqn:lem-greedy3-3}
\end{align}
\RED{%
By the choice of $(\sigma,X)$,
where 
$\sigma =+1$ in our first case, 
}% 
we have
$g(p+ \unitvec{X})  \leq g(p + \unitvec{X \setminus B})$,
which, together with (\ref{eqn:lem-greedy3-1})--(\ref{eqn:lem-greedy3-3}),  
implies that
$g(q^*)  \geq g(q^* + \unitvec{B \cap X})$,
i.e., $q^* + \unitvec{B \cap X} \in \arg\min g$.

Since $\min_{i \in N}\{q^*(i) - p(i)\} \leq 0 <
\max_{i \in N}\{q^*(i) - p(i)\}$
by the assumption and (\ref{eqn:claim1})
\RED{%
established in Step 2,
}%,
we have $A \cap B = \emptyset$, which, together with 
(\ref{eqn:claim2})
\RED{%
established in Step 3,
}%,
implies $A \subseteq X \setminus B$.
 Hence, it holds that
\begin{equation*}
 \|(q^* + \unitvec{B \cap X}) - (p + \unitvec{X})  \|^+_\infty
= \|q^* - p - \unitvec{X\setminus B} \|^+_\infty
= \|q^* - p  \|^+_\infty - 1.
\end{equation*}
 We also have
\begin{equation}
\|(q^* + \unitvec{B \cap X}) - (p + \unitvec{X})  \|^-_\infty
=
 \|q^*  - p - \unitvec{X\setminus B}  \|^-_\infty
= \|q^* - p  \|^-_\infty,
\end{equation}
where the second equality follows from the definition of $B$.
 Hence, it holds that
\begin{align*}
\mu(p + \unitvec{X})
 & \leq 
 \|(q^* + \unitvec{B \cap X}) - (p + \unitvec{X})  \|^+_\infty
 + \|(q^* + \unitvec{B \cap X}) - (p + \unitvec{X})  \|^-_\infty
 \nonumber \\
 & =
 ( \|q^* - p  \|^+_\infty - 1) + \|q^* - p  \|^-_\infty
 \ = \ \mu(p) -1.
\end{align*}
\RED{%
We have completed the proof of  Lemma \ref{lem:Lnat-greedy}
when $\sigma =+1$.
}%

\medskip

\RED{%
Next we go on to the second case with $\sigma =-1$.
}%
It follows from the definition of an s.s.~quasi \Lnat-convex function
that a function $g$ is s.s.~quasi \Lnat-convex
 if and only if the function $h(p):=g(-p)$ is s.s.~quasi \Lnat-convex.
This symmetry implies that 
the proof for the case of $\sigma = - 1$ 
can be done symmetrically 
to the proof for the case of $\sigma = + 1$.
The following paragraphs for the case of $\sigma = - 1$ 
are precisely this symmetric argument,
and do not contain anything essentially different 
from the case of $\sigma = + 1$.

 {\bf [Step 1 ($\sigma = - 1$)]} \quad
 We first show the inequality
$\mu(p - \unitvec{X}) \geq \mu(p) -1$.
 For every $d \in \Z^n$ and $Y \subseteq N$, we have
\begin{equation*}
    \|d + \unitvec{Y} \|^+_\infty \geq  
\RED{%
\|d  \|^+_\infty,  
}%
\qquad
    \| d + \unitvec{Y} \|^-_\infty \geq 
\RED{%
 \| d  \|^-_\infty -1.
}%
\end{equation*}
 Hence, it holds that
\begin{eqnarray*}
\mu(p - \unitvec{X})
& = &
\min\{ \|q - (p - \unitvec{X}) \|^+_\infty
+ \|q - (p - \unitvec{X}) \|^-_\infty
\mid q \in \arg\min g \}
\nonumber \\
& \geq &
\min\{ \| q - p \|^+_\infty +  \| q - p \|^-_\infty
 \mid q \in \arg\min g \} - 1
\nonumber \\
& = &  \mu(p) -1.
\end{eqnarray*}

\RED{%
In the following (i.e., in Steps 2 to 4),
}%
we prove the reverse inequality
\begin{equation}
\label{eqn:rev-ineqM} 
\mu(p - \unitvec{X}) \leq \mu(p) -1.
\end{equation} 
The outline of the proof is as follows.
 We denote
\begin{align}
  S & = \{ q \in \arg\min g \mid
 \|q - p \|^+_\infty +  \| q - p \|^-_\infty = \mu(p) \},
\\
 \xi & = \max\{\|q - p\|_\infty^- \mid q \in S \}.
\end{align}
 Let  $q^*$ be a vector in $S$
with $\|q^* - p \|^-_\infty = \xi$,
and assume that $q^*$ is 
\RED{%
a maximal vector
}%
among all such vectors.
In Step 2,
we show that 
\RED{%
\begin{equation}
\label{eqn:claim1M} 
\xi > 0,  \quad  \mbox{i.e.,} \quad  \|q^* - p \|^-_\infty > 0.
\end{equation}
}%
 Note that this condition is equivalent to
$\suppm(q^* - p)\neq \emptyset$.
 Using this, we then prove,
in Step 3, that
\RED{%
\begin{equation}
\label{eqn:claim2M}
\arg\max_{i \in N}\{p(i)- q^*(i) \} \subseteq X. 
\end{equation}
}%
 By using (\ref{eqn:claim1M}) and (\ref{eqn:claim2M}), 
we derive the inequality (\ref{eqn:rev-ineqM})
in Step 4.

 {\bf [Step 2  ($\sigma = - 1$): proof of (\ref{eqn:claim1M})]} \quad
 Assume, to the contrary, that 
\RED{%
$(\xi=) \ \|q^* - p \|^-_\infty = 0$, 
}%
i.e.,
\RED{%
$q^* \geq  p$ 
}%
holds. 
 This assumption implies $\|q^* - p \|^+_\infty > 0$
since $\mu(p) > 0$.
\BLU{%
By Lemma \ref{lem:Lnat-subm}
for $(q^*, p - \unitvec{X})$,
}%%
one of the following three conditions holds:
\BLU{%
\begin{align}
& 
 g(p - \unitvec{X}) > g(((p - \unitvec{X}) \wedge q^*) + \unitvec{Z}),
 \label{eqn:lem-greedy7-1-1M}
\\
& 
 g(q^*) > g(((p - \unitvec{X}) \vee q^*) - \unitvec{Z}),
 \label{eqn:lem-greedy7-1-2M}
\\
& 
 g(p - \unitvec{X}) = g(((p - \unitvec{X}) \wedge q^*) + \unitvec{Z}) 
\mbox{ and } g(q^*) = g(((p - \unitvec{X}) \vee q^*) - \unitvec{Z}),
 \label{eqn:lem-greedy7-1-3M}
\end{align}
}%
where $Z= \suppp( q^* - (p - \unitvec{X}))$.
Let 
$Y= \{i \in N \mid  q^*(i) - p(i) = 0\}$, 
which may be the empty set.
 Since 
\BLU{%
$q^* \geq p$ by the assumption
and $Z = (N \setminus Y) \cup X$,
}%% 
\BLU{%
 we have
\begin{align*}
& 
((p - \unitvec{X}) \wedge q^*) + \unitvec{Z}
 = (p - \unitvec{X}) + \unitvec{Z}
 = p + \unitvec{N \setminus (X\cup Y)},
\\ &
((p - \unitvec{X}) \vee q^*) - \unitvec{Z}=
 q^* - \unitvec{Z}=
q^* - \unitvec{(N \setminus Y) \cup X}.
\end{align*}
}%color
With these equations, the three conditions
 (\ref{eqn:lem-greedy7-1-1M})--(\ref{eqn:lem-greedy7-1-3M})
can be rewritten as follows:
\begin{align}
& 
  g(p - \unitvec{X}) >  g(p + \unitvec{N \setminus (X \cup Y)}),
 \label{eqn:lem-greedy7-2-1M}
\\
& 
  g(q^*)  >   g(q^* - \unitvec{(N \setminus Y) \cup X}),
 \label{eqn:lem-greedy7-2-2M}
\\
& 
  g(p - \unitvec{X}) =  g(p + \unitvec{N \setminus (X \cup Y)})
\mbox{ and }
  g(q^*)  =   g(q^* - \unitvec{(N \setminus Y) \cup X}).
 \label{eqn:lem-greedy7-2-3M}
\end{align} 
By the choice of $(\sigma,X)$,
where 
$\sigma =-1$ in our second case, 
we have
$g(p - \unitvec{X}) \leq
g(p + \unitvec{N \setminus (X \cup Y)})$.
 From this and (\ref{eqn:lem-greedy7-2-1M})--(\ref{eqn:lem-greedy7-2-3M})
 follows that
$ g(q^*)  \geq
g(q^* - \unitvec{(N \setminus Y) \cup X})$,
implying that
$q^* - \unitvec{(N \setminus Y) \cup X} \in \arg\min g$.
 By $q^* \geq p$ and the definition of $Y$, we have
\begin{equation}
 \label{eqn:lem-greedy7-11M} 
\|(q^* - \unitvec{(N \setminus Y) \cup X}) - p\|_\infty^+
 =  
\max_{i \in N \setminus Y}
\RED{%
\{ (q^*(i)-1) - p(i)\}
}%
 = \|q^* - p\|_\infty^+ -1
\end{equation}
since $\|q^* - p\|_\infty^+ > 0$.
 We also have
\begin{align}
\|(q^* - \unitvec{(N \setminus Y) \cup X}) - p\|_\infty^-
& \leq  1\  =\  \|q^* - p\|_\infty^- +1 
 \label{eqn:lem-greedy7-3M} 
\end{align}
since $\|q^* - p\|_\infty^- = 0$.
 From (\ref{eqn:lem-greedy7-11M}) and (\ref{eqn:lem-greedy7-3M}) follows that
\begin{align}
\mu(p)
& \leq 
\|(q^* - \unitvec{(N \setminus Y) \cup X})  - p \|^+_\infty
 + \|(q^* - \unitvec{(N \setminus Y) \cup X}) - p\|_\infty^-
\notag\\
&\leq   \|q^*  - p \|^+_\infty + \|q^* - p\|_\infty^-
 =  \mu(p),
 \label{eqn:lem-greedy7-12M} 
\end{align}
where the first inequality is by 
the definition of $\mu(p)$.
 Hence, the inequality (\ref{eqn:lem-greedy7-3M}) 
and the first inequality in (\ref{eqn:lem-greedy7-12M}) must hold with equality,
i.e., we have $q^* - \unitvec{(N \setminus Y) \cup X} \in S$ 
and 
\begin{equation*}
 \|(q^* - \unitvec{(N \setminus Y) \cup X}) - p\|_\infty^-
 =  \|q^* - p\|_\infty^- +1 
>  \|q^* - p\|_\infty^- = \xi. 
\end{equation*}
 This, however, is a contradiction to the definition of $\xi$.
 Hence, (\ref{eqn:claim1M}) holds.

 {\bf [Step 3  ($\sigma = - 1$): proof of (\ref{eqn:claim2M})]} \quad
 We denote
\begin{equation*}
A = \arg\max_{i \in N}
\RED{%
\{ p(i) - q^*(i) \}.
}%
\end{equation*}
 Then, (\ref{eqn:claim2M}) is simply rewritten as $A \subseteq X$.

 Assume, to the contrary, that $A \setminus X \neq \emptyset$ holds.
\RED{%
 We  will show 
}%
that
$q^* + \unitvec{A \setminus X} \in \arg\min g$.
 By (\ref{eqn:claim1M})
\RED{%
established in Step 2,
}%
it holds that  $\xi=\|q^* - p\|_\infty^-   > 0$.
 Therefore, we have $A \subseteq \suppm(q^* - p)$, from which follows that 
\begin{equation*}
 \suppm(q^* - (p - \unitvec{X})) \supseteq A \setminus X
\neq \emptyset.
\end{equation*}
 Since $A \setminus X \neq \emptyset$, we also have
\begin{equation*}
 \arg\max_{i \in N}
\RED{%
\{(p - \unitvec{X})(i) - q^*(i) \} 
}%
= A \setminus X.
\end{equation*}
 Hence, Lemma \ref{lem:Lnat-exchange} 
for $(p - \unitvec{X}, q^*)$,
together with the relation
$(p - \unitvec{X}) - \unitvec{A \setminus X} = p - \unitvec{X \cup A}$,
implies that
one of the following three conditions holds:
\begin{align}
&
   g(q^*)  > g(q^* + \unitvec{A \setminus X}),
 \label{eqn:lem-greedy2-1M}
\\
&
 g(p - \unitvec{X}) > g(p - \unitvec{X \cup A}),
 \label{eqn:lem-greedy2-2M}
\\
&
   g(q^*)  = g(q^* + \unitvec{A \setminus X})
\mbox{ and }
 g(p - \unitvec{X}) = g(p - \unitvec{X \cup A}).
 \label{eqn:lem-greedy2-3M}
\end{align}
 By the choice of $(\sigma,X)$,
where 
$\sigma =-1$ in our second case, 
 we have
$g(p - \unitvec{X}) \leq  g(p - \unitvec{X \cup A})$,
which, together with (\ref{eqn:lem-greedy2-1M})--(\ref{eqn:lem-greedy2-3M}), 
implies that
 $g(q^*) \geq g(q^* + \unitvec{A \setminus X})$, i.e.,
$q^* + \unitvec{A \setminus X} \in \arg\min g$.

 Since $A  \setminus X \subseteq A \subseteq \suppm(q^* - p)$, we have
\begin{align}
 \|(q^* + \unitvec{A \setminus X}) - p\|_\infty^-
& \leq  \|q^* - p\|_\infty^- = \xi,
 \label{eqn:lem-greedy7-5M} 
\\
\|(q^* + \unitvec{A \setminus X}) - p\|_\infty^+
& = \|q^* - p\|_\infty^+,
%%\notag
% \label{eqn:lem-greedy7-8M} 
\end{align}
 from which follows that
\begin{align}
 \mu(p)
& \leq \|(q^* + \unitvec{A \setminus X}) - p\|_\infty^+
 +
 \|(q^* + \unitvec{A \setminus X}) - p\|_\infty^- \notag\\
& \leq \|q^* - p\|_\infty^+
 + \|q^* - p\|_\infty^- = \mu(p),
  \label{eqn:lem-greedy7-6M} 
\end{align}
where the first inequality is by the definition of $\mu(p)$. 
 Hence, the inequality (\ref{eqn:lem-greedy7-5M})
and the first inequality in (\ref{eqn:lem-greedy7-6M}) 
must hold with equality.
 Hence, the vector $q^* + \unitvec{A \setminus X}$
belongs to $S$ with
$\|(q^* + \unitvec{A \setminus X}) - p\|_\infty^-  = \xi$,
a contradiction to the 
\RED{%
maximality of $q^*$.
}%

 {\bf [Step 4  ($\sigma = - 1$): proof of (\ref{eqn:rev-ineqM})]} \quad
 To derive the inequality (\ref{eqn:rev-ineqM})
from (\ref{eqn:claim1M}) and (\ref{eqn:claim2M}),
 we distinguish two cases:
$\min_{i \in N}
\RED{%
\{ p(i) - q^*(i) \}
}%
 > 0$ 
or
$\min_{i \in N}
\RED{%
\{ p(i) - q^*(i) \}
}%
 \leq 0$.

We first consider
the case with 
$\min_{i \in N}
\RED{%
\{ p(i) - q^*(i) \}
}%
 > 0$.
%% or $B \cap X = \emptyset$ holds.
 Since $q^*(i) < p(i)$ for all $i \in N$, it holds that
$q^* \leq p - \unitvec{X}$.
 Therefore, we have
\begin{equation*}
  \|q^* - (p - \unitvec{X}) \|^+_\infty
=0= \|q^* - p \|^+_\infty.
\end{equation*}
 By (\ref{eqn:claim2M})
\RED{%
established in Step 3,
}%,
it holds that
\begin{equation*}
  \|q^* - (p - \unitvec{X}) \|^-_\infty
= \|q^* - p \|^-_\infty - 1.
\end{equation*}
 Therefore, it follows that
\begin{align*}
\mu(p - \unitvec{X}) 
& \leq 
  \|q^* - (p - \unitvec{X}) \|^+_\infty
+ \|q^* - (p - \unitvec{X}) \|^-_\infty
\nonumber
\\
& = 
\RED{%  
\|q^* - p  \|^+_\infty +  (\|q^* - p  \|^-_\infty -1) 
} %
 =  \mu(p) -1.
\end{align*}

We next consider the remaining case where
 $\min_{i \in N}
\RED{%
\{ p(i) - q^*(i) \}
}%
\leq 0$.
%% and $B \cap X \neq \emptyset$.
 We denote
\begin{equation*}
 B = \arg\min_{i \in N}
\RED{%
\{ p(i) - q^*(i) \}.
}%
\end{equation*}
\RED{%
 We  will show 
}%
that
$q^* - \unitvec{B \cap X} \in \arg\min g$ holds.
 If $B \cap X = \emptyset$, then
$q^* - \unitvec{B \cap X} = q^*\in \arg\min g$.
 Hence, we assume
$B \cap X \neq \emptyset$.
 Since
\begin{equation*}
\max_{i \in N}
\RED{%
\{ q^*(i) - p(i) \} 
}%
\geq 0,
\qquad  B = \arg\max_{i \in N}
\RED{%
\{ q^*(i) - p(i) \} ,
}%
\end{equation*}
it holds that
\begin{equation}
\RED{%
\suppp( q^* - (p - \unitvec{X})) 
}%
\supseteq B \cap X \neq \emptyset, 
\quad
 \arg\max_{i \in N} 
\RED{%
\{ q^*(i) - (p - \unitvec{X})(i) \} 
}%
= B \cap X.
\end{equation}
 It follows from Lemma \ref{lem:Lnat-exchange} 
for $(q^*, p - \unitvec{X})$
and the relation
$(p - \unitvec{X}) + \unitvec{B \cap X} = p - \unitvec{X \setminus B}$
that one of the following three conditions holds:
\begin{align}
&
   g(p - \unitvec{X})    >  g(p - \unitvec{X \setminus B}),
 \label{eqn:lem-greedy3-1M}
\\
&
g(q^*)  > g(q^* - \unitvec{B \cap X}),
 \label{eqn:lem-greedy3-2M}
\\
&
   g(p - \unitvec{X})    =  g(p - \unitvec{X \setminus B}) 
\mbox{ and }
g(q^*)  = g(q^* - \unitvec{B \cap X}).
 \label{eqn:lem-greedy3-3M}
\end{align}
By the choice of $(\sigma,X)$,
where 
$\sigma =-1$ in our second case, 
we have
$g(p - \unitvec{X})  \leq g(p - \unitvec{X \setminus B})$,
which, together with (\ref{eqn:lem-greedy3-1M})--(\ref{eqn:lem-greedy3-3M}),  
implies that
$g(q^*)  \geq g(q^* - \unitvec{B \cap X})$,
i.e., $q^* - \unitvec{B \cap X} \in \arg\min g$.

Since 
\RED{%
$\min_{i \in N} \{ p(i) - q^*(i) \}
\leq 0 < \max_{i \in N} \{ p(i) - q^*(i) \}$
}%
by the assumption and (\ref{eqn:claim1M}) 
\RED{%
established in Step 2,
}%
we have $A \cap B = \emptyset$, which, together with 
(\ref{eqn:claim2M})
\RED{%
established in Step 3,
}%
implies $A \subseteq X \setminus B$.
 Hence, it holds that
\begin{equation*}
 \|(q^* - \unitvec{B \cap X}) - (p - \unitvec{X})  \|^-_\infty
= \|q^* - p + \unitvec{X\setminus B} \|^-_\infty
= \|q^* - p  \|^-_\infty - 1.
\end{equation*}
 We also have
\begin{equation*}
\|(q^* - \unitvec{B \cap X}) - (p - \unitvec{X})  \|^+_\infty
=
 \|q^*  - p + \unitvec{X\setminus B}  \|^+_\infty
= \|q^* - p  \|^+_\infty,
\end{equation*}
where the second equality follows from the definition of $B$.
 Hence, it holds that
\begin{align*}
\mu(p - \unitvec{X})
 & \leq 
 \|(q^* - \unitvec{B \cap X}) - (p - \unitvec{X})  \|^+_\infty
 + \|(q^* - \unitvec{B \cap X}) - (p - \unitvec{X})  \|^-_\infty
 \nonumber \\
 & 
\RED{\, 
= \|q^* - p  \|^+_\infty +  ( \|q^* - p  \|^-_\infty - 1) 
}%
 \ = \ \mu(p) -1.
\end{align*}

We have completed the proof of  Lemma \ref{lem:Lnat-greedy} when 
\BLU{%
$\sigma = - 1$,
}%
and hence the proof of  Lemma \ref{lem:Lnat-greedy}
for both cases with $\sigma = \pm 1$.
Thus Theorem \ref{thm:greedy-iter} is proved.

\subsection{Proof of Theorem \ref{thm:ascend-iter}}
\label{SCprfTHascend-iter}

 The proof of Theorem \ref{thm:ascend-iter} 
is quite similar to and simpler than
that of Theorem \ref{thm:greedy-iter}.
 Theorem \ref{thm:ascend-iter} can be proved by
using the following property repeatedly.

\begin{lemma}
\label{lem:ascend-iteration}
 Let $p \in \Z^n$ be a vector with $\hat\mu(p) > 0$,
and $X \subseteq N$ be a set that minimizes
the value of $g(p+ \unitvec{X})$.
 Then, $\hat\mu(p+  \unitvec{X})= \hat\mu(p)-1$.
\end{lemma}

\RED{%
In the rest of this section,
}%
we give a proof of Lemma \ref{lem:ascend-iteration}.
 The inequality $\hat{\mu}(p + \unitvec{X}) \geq \hat{\mu}(p) -1$ can be
 shown as follows. 
 By the triangle inequality,
we have $\|q - (p + \unitvec{X})  \|_\infty \geq \|q - p  \|_\infty - 1$
for every $q \in \Z^n$.
 Taking the minimum over all $q \in \arg\min g$
with $q \geq p + \unitvec{X}$, we obtain
\begin{align*}
 \hat{\mu}(p + \unitvec{X})
 & \geq
 \min\{ \|q - p  \|_\infty
 \mid q \in \arg\min g,\ q \geq p + \unitvec{X} \} - 1\\
 & \geq 
 \min\{ \|q - p  \|_\infty
 \mid q \in \arg\min g,\ q \geq p \} - 1
  =   \hat{\mu}(p) -1.
\end{align*}
 In the following, we show the reverse inequality:
\begin{equation}
\label{eqn:rev-ineq2}
 \hat\mu(p + \unitvec{X}) \leq \hat\mu(p)-1. 
\end{equation}

 Let $p^*$ be a vector such that
$p^* \in \arg\min g$, $p^* \geq p$, and
$\|p^* - p\|_\infty = \hat{\mu}(p)$,
and assume that $p^*$ is minimal among all such vectors.
 We denote
\[
A = \arg\max_{i \in N}\{p^*(i) - p(i)\}.
\]
 We have $p^* \neq p$ and $\max_{i \in N}\{p^*(i) - p(i)\}>0$
since $\|p^* - p\|_\infty = \hat{\mu}(p) > 0$ and $p^* \geq p$.

 We claim that
\begin{equation}
\label{eqn:lem-greedy4L} 
A \subseteq X.
\end{equation}
 Assume, to the contrary, that $A \setminus X \neq \emptyset$ holds.
 Since $A \subseteq \suppp(p^* - p)$,
we have
\[
 \suppp(p^* - (p + \unitvec{X})) \supseteq A \setminus X
\neq \emptyset.
\]
 We also have
\[
 \arg\max_{i \in N}\{p^*(i) - (p + \unitvec{X})(i)\} = A \setminus X.
\]
 Hence, Lemma \ref{lem:Lnat-exchange} 
\RED{%
for $(p^*, p+ \unitvec{X})$,
together with the relation
$(p + \unitvec{X}) + \unitvec{A \setminus X} = p + \unitvec{X \cup A}$,
implies that
one of the following three conditions holds:
}%
\begin{align}
 &
   g(p^*)  > g(p^* - \unitvec{A \setminus X}),
 \label{eqn:lem-greedy2L-1}
 \\
 &
 g(p+ \unitvec{X})  >  g(p + \unitvec{X \cup A}),
 \label{eqn:lem-greedy2L-2}
 \\
 &
   g(p^*)  = g(p^* - \unitvec{A \setminus X})
 \mbox{ and }
 g(p+ \unitvec{X})  =  g(p + \unitvec{X \cup A}).
 \label{eqn:lem-greedy2L-3}
\end{align}
 By the choice of $X$, we have
$g(p+ \unitvec{X}) \leq  g(p + \unitvec{X \cup A})$.
 This inequality, together with (\ref{eqn:lem-greedy2L-1})--(\ref{eqn:lem-greedy2L-3}),  
implies that
$g(p^*) \geq g(p^* - \unitvec{A \setminus X})$,
i.e., $p^* - \unitvec{A \setminus X} \in \arg\min g$ holds.
 This, however, is a contradiction to the choice of $p^*$
since 
\[
 p^* \geq p^* - \unitvec{A \setminus X} \geq p,
\qquad
\|(p^* - \unitvec{A \setminus X}) - p\|_\infty \leq
\|p^*  - p\|_\infty = \hat{\mu}(p).
\]
 Hence, we have (\ref{eqn:lem-greedy4L}).

 We now prove the inequality (\ref{eqn:rev-ineq2}). 
\RED{%
We distinguish two cases, 
depending on $p^* \geq p + \unitvec{X}$ or not.
If $p^* \geq p + \unitvec{X}$ is true, 
}%
we have 
\[
 \hat{\mu}(p + \unitvec{X})
\leq \|p^* - (p + \unitvec{X})\|_\infty
= \|p^* - p\|_\infty - 1
= \hat{\mu}(p) -1,
\]
where the first equality is by (\ref{eqn:lem-greedy4L}).

\RED{%
If the condition $p^* \geq p + \unitvec{X}$ fails,
we have
}%
$B \cap X \neq \emptyset$ for
$B = \{i \in N \mid p^*(i) = p(i)\}$.
 Since $\max_{i \in N}\{p^*(i) - p(i)\} > 0$,
we have $A \cap B = \emptyset$, which, together with
(\ref{eqn:lem-greedy4L}), implies 
$A \subseteq X \setminus B$.
 Since $p^* \geq p$, we have
\begin{equation}
\label{eqn:Lgreedy1}
  p^*(i) = p(i) \quad (\forall i \in B),
\qquad
 p^*(i) > p(i) \quad (\forall i \in N \setminus B),
\end{equation}
from which $p^* + \unitvec{B \cap X} \geq p+ \unitvec{X}$ follows.
 As shown below, we have
$p^* + \unitvec{B \cap X} \in \arg\min g$.
 Hence, it holds that
\begin{align*}
   \hat{\mu}(p + \unitvec{X})
 \leq  \|(p^* + \unitvec{B \cap X}) - (p + \unitvec{X})\|_\infty
 & =  \|p^* - p - \unitvec{X \setminus B} \|_\infty\\
 & =  \|p^* - p\|_\infty - 1
  =  \hat{\mu}(p) -1,
\end{align*}
where the second equality is 
by $A \subseteq X \setminus B$, (\ref{eqn:Lgreedy1}),
and the definition of $A$.

 We now show that $p^* + \unitvec{B \cap X} \in \arg\min g$ holds.
 The condition (\ref{eqn:Lgreedy1}) implies 
\[
 \suppp((p+ \unitvec{X})- p^*) = 
 \arg\max_{i \in N}\{(p+ \unitvec{X})(i)- p^*(i)\} = B \cap X.
\]
 Hence, it follows from Lemma \ref{lem:Lnat-exchange} 
\RED{%
for $(p+ \unitvec{X}, p^*$),
as well as the relation
$(p + \unitvec{X}) - \unitvec{B \cap X} = p + \unitvec{X \setminus B}$,
that one of the following three conditions holds:
}%
\begin{align}
&
    g(p+ \unitvec{X}) >  g(p + \unitvec{X \setminus B}), 
 \label{eqn:lem-greedy3L-1}
\\
&
g(p^*)  > g(p^* + \unitvec{B \cap X}),
 \label{eqn:lem-greedy3L-2}
\\
&
    g(p+ \unitvec{X})   =  g(p + \unitvec{X \setminus B}) 
\mbox{ and }
g(p^*)  = g(p^* + \unitvec{B \cap X}).
 \label{eqn:lem-greedy3L-3}
\end{align}
 By the choice of $X$, we have
$g(p+ \unitvec{X})  \leq g(p + \unitvec{X \setminus B})$,
which, together with (\ref{eqn:lem-greedy3L-1})--(\ref{eqn:lem-greedy3L-3}),  
implies that
$g(p^*)  \geq g(p^* + \unitvec{B \cap X})$,
i.e., $p^* + \unitvec{B \cap X}$ is a minimizer of $g$.

 This concludes the proof of Lemma \ref{lem:ascend-iteration}
(and also of Theorem \ref{thm:ascend-iter}).


\begin{thebibliography}{99}


\bibitem{Mdcasiam} 
Murota, K.:
Discrete Convex Analysis.
Society for Industrial and Applied Mathematics, Philadelphia (2003)




\bibitem{MS03quasi} 
Murota, K., Shioura, A.:
Quasi M-convex and L-convex functions:
quasi-convexity in discrete optimization.
Discrete Applied Mathematics {\bf 131} 467--494 (2003)



\bibitem{MS14exbndLmin}
Murota, K., Shioura, A.:
Exact bounds for steepest descent algorithms of L-convex function minimization.
Operations Research Letters {\bf 42}, 361--366 (2014) 
\end{thebibliography}
\end{document}